\documentclass[12pt,a4paper]{amsart}
\usepackage[a4paper,textwidth=16cm,textheight=23cm]{geometry}
\usepackage[utf8]{inputenc}
\usepackage{amsmath}
\usepackage{amsfonts}
\usepackage{amssymb}
\usepackage{xcolor}
\usepackage{graphicx}
\providecommand{\U}[1]{\protect\rule{.1in}{.1in}}
\newtheorem{theorem}{Theorem}[section]

\newtheorem{corollary}[theorem]{Corollary}

\newtheorem{example}[theorem]{Example}

\newtheorem{lemma}[theorem]{Lemma}

\newtheorem{proposition}[theorem]{Proposition}

\graphicspath{{Slike/}}

\newcommand{\F}{\mathbb F}
\newcommand{\Mat}{M_2(\F_q)}
\newcommand{\GL}{\mathrm{GL}}
\newcommand{\D}{D_q}
\newcommand{\rank}{\mathrm{rank}}
\DeclareMathOperator{\PG}{PG}

\title[Domination versus covering in total graphs of finite rings]
  {Domination versus covering in total graphs of finite rings}

\author{David Dol\v{z}an, Riste \v{S}krekovski}
\address{D.~Dol\v{z}an: Department of Mathematics, Faculty of Mathematics and Physics,
University of Ljubljana, Jadranska 19, SI-1000 Ljubljana, Slovenia,
and Institute of Mathematics, Physics and Mechanics,
Jadranska 19, SI-1000 Ljubljana, Slovenia}
\email{david.dolzan@fmf.uni-lj.si}

\address{R.~\v{S}krekovski: Department of Mathematics, Faculty of Mathematics and Physics,
University of Ljubljana, Jadranska 19, SI-1000 Ljubljana, Slovenia,
and Faculty of Information Studies, University of Novo Mesto, Ljubljanska cesta 31a, SI-8000 Novo mesto, Slovenia}
\email{riste.skrekovski@fmf.uni-lj.si}


\thanks{D.~Dol\v{z}an acknowledges partial support from the
Slovenian Research and Innovation Agency program P1-0222. R.~\v{S}krekovski acknowledges partial support from the Slovenian Research and Innovation Agency program P1-0383.}

\begin{document}

\begin{abstract}
For a finite ring $R$, let $\kappa(R)$ be the minimum number of
translates of its nonunits covering $R$. We prove that $\kappa$ is
invariant modulo the Jacobson radical and takes the minimum over
direct factors. It equals the domination number of the total graph
when the radical is nonzero and for $M_n(\F_q)$ with $n\ge3$.
For odd $q$, the domination number of $M_2(\F_q)$ is
$\min\{\kappa(M_2(\F_q)),2q-2\}$, disproving the
Dol\v{z}an--Oblak conjecture [\emph{Comm. Algebra} \textbf{43} (2015) 2903--2911]. We obtain geometric lower bounds and
determine the covering, domination, and total domination numbers
of $M_2(\F_q)$ for $q=3,4,5,7$.
\end{abstract}

\maketitle

\section{Introduction}
The total graph of a ring combines an additive relation with a multiplicative obstruction: two distinct ring elements are adjacent when their sum is a zero-divisor. For a full matrix ring over a finite field, this means that two distinct matrices are adjacent precisely when their sum is singular. Questions about domination can therefore be approached through determinants and the geometry of the singular matrices.

The total graph was introduced for commutative rings by Anderson and Badawi \cite{AB}. Throughout, $R$ denotes a finite nonzero ring with identity, and
$Z(R)$ denotes its set of zero-divisors, including zero. In a
finite ring, the left and right zero-divisors coincide and are
precisely the nonunits. The \emph{total graph} $\tau(R)$ is the simple graph with vertex set $R$ in which distinct vertices $x,y$ are adjacent if and only if $x+y\in Z(R)$. A subset $D\subseteq V(G)$ is \emph{dominating} if every vertex of $V(G)\setminus D$ has a neighbour in $D$. The minimum cardinality of such a set is the \emph{domination number} $\gamma(G)$. In particular, the selected vertices are not required to have neighbours within the dominating set.

Domination in the total graph of a finite commutative ring was studied by Shekarriz, Shirdareh Haghighi, and Sharif \cite{SH}. Dol\v{z}an and Oblak \cite{DO} considered total graphs of arbitrary finite rings, including noncommutative rings, and investigated connectivity, Hamiltonicity, and domination. Full matrix rings are natural in this setting because the semisimple quotient of a finite ring is a direct product of matrix rings over finite fields.

One construction in \cite[Lemma~4.2]{DO} gives the dominating set
\[
 \{0\}\cup\{xE_{1j}:x\in\F_q^\times,\ 1\leq j\leq n\}
\]
in $\tau(M_n(\F_q))$ for $n\geq2$. Here $E_{ij}$ denotes the usual matrix unit. This set has cardinality $n(q-1)+1$. The underlying argument uses the linear dependence of the determinant on a single matrix entry: changing an appropriate first-row entry makes an invertible matrix singular. This gives a natural upper bound whose possible sharpness motivates the conjecture considered below.

The conjecture in \cite[p.~2910]{DO} proposes that, if a finite ring $R$ satisfies
\[
 R/J(R)\cong\prod_{i=1}^k M_{n_i}(\F_{q_i}),\qquad n_i\geq2,
\]
then
\begin{equation}
 \gamma(\tau(R))=\min_{1\leq i\leq k}\{n_i(q_i-1)+1\}.
\label{For_conjecture}
\end{equation}
Here $J(R)$ is the Jacobson radical. In particular, this predicts $\gamma(\tau(M_2(\F_q)))=2q-1$.

We show that this prediction fails for every odd prime power $q$. Our construction uses the affine line $\{I+tE_{12}:t\in\F_q\}$ together with suitable scalar matrices, where $I$ is the identity matrix. The resulting dominating set has $2q-2$ elements. The set is dominating because it contains every matrix whose sum with each selected matrix is invertible. A complementary lower bound follows by fixing the first row
and counting the possible second rows.

\begin{theorem}
\label{Tm_main}Let $q$ be a prime power.
\begin{enumerate}
\item[(a)] For every prime power $q\ge3$,
\[
\gamma(\tau(M_2(\F_q)))\ge q+1.
\]
If $q$ is odd, then $\gamma(\tau(M_2(\F_q)))\le 2q-2$.
Consequently, (\ref{For_conjecture}) fails for every ring $M_2(\F_q)$ of odd field order.
\item[(b)] If $q$ is odd, then
\[
\gamma(\tau(\Mat))=\min\{\kappa(\Mat),2q-2\},
\]
where $\kappa$ is the covering number defined in Section~\ref{Sec_covering}.
\item[(c)] For $q=3,4,5,7$, the domination numbers of $\tau(\Mat)$ are
$4,7,8,12$, respectively, while
$\kappa(\Mat)=\gamma_t(\tau(\Mat))=2q-1$ in each case.
\end{enumerate}
\end{theorem}

Part~(a) is proved in Section~\ref{Sec_domination}. Parts~(b)
and~(c) are proved at the end of Section~\ref{Sec_fibres}.

The gap between domination and covering in Theorem~\ref{Tm_main}(b)
is a phenomenon of $2\times2$ matrices.

\begin{theorem}\label{Tm_n3_intro}
Let $n\ge3$ and let $q$ be a prime power. Then every dominating set of
$\tau(M_n(\F_q))$ with at most $n(q-1)$ elements is a singular-sum
covering set, and
\[
\gamma(\tau(M_n(\F_q)))=\kappa(M_n(\F_q)).
\]
Consequently, for $n\ge3$, the value predicted by~\eqref{For_conjecture}
for $M_n(\F_q)$ is correct if and only if
$\kappa(M_n(\F_q))=n(q-1)+1$.
\end{theorem}

Theorem~\ref{Tm_n3_intro} is proved in Section~\ref{Sec_n3}. There we
also extend it to every finite ring $R$ as in~\eqref{For_conjecture}
other than $M_2(\F_q)$: for all such rings
$\gamma(\tau(R))=\kappa(R)=\min_i\kappa(M_{n_i}(\F_{q_i}))$
(Corollary~\ref{Cor_products}). Thus, among the rings in~\eqref{For_conjecture}, domination and
covering can differ only for $M_2(\F_q)$ with $q$ odd.
For $q=3$, the two bounds in Theorem~\ref{Tm_main}(a) coincide. For $q\ge7$, our
geometric estimates give the stronger bound
\[
\gamma(\tau(\Mat))\ge q+\lceil\sqrt q\rceil+2
\]
by Corollary~\ref{cor:sqrt}, Theorem~\ref{Tm_main}(b), and
Proposition~\ref{Prop_three}.

Total domination imposes an additional requirement: each selected
vertex must itself have a distinct selected neighbour. We study this
parameter through a covering condition that requires a zero-divisor sum
at every ring element, while allowing a selected element to cover itself
when its double is a zero-divisor. This intermediate condition agrees
with domination in characteristic two. For nonlocal rings of odd order, we show that its minimum
size equals the total domination number. The same equality holds in every characteristic for a direct
product of two finite nonzero rings. We also prove that the domination number equals the covering number whenever the Jacobson radical is nonzero. Consequently,
a strict gap between domination and covering
can occur only for semisimple rings. We denote this covering number by $\kappa(R)$.
For $M_n(\F_q)$, this covering number is the minimum size of a
rank-metric covering code of radius at most $n-1$~\cite{GY08,GY09}.
In particular, Theorem~\ref{Tm_main}(c) gives domination number $12$
and covering and total domination numbers $13$ for $M_2(\F_7)$.

Section~\ref{Sec_domination} establishes the counterexamples, and
Section~\ref{Sec_covering} proves the radical and product formulas.
Section~\ref{Sec_geometric} bounds the covering number by means of
planar blocking sets. Section~\ref{Sec_fibres} relates this bound to
 domination by counting uncovered matrices in affine planes.
Section~\ref{Sec_n3} shows that domination and covering coincide for
matrices of size at least three.

\section{Counterexamples and bounds for domination}\label{Sec_domination}

\begin{proposition}
\label{Pr_upper}Let $q$ be an odd prime power. The set
\begin{equation}
 \D=\{0\}\cup\{I+tE_{12}:t\in\F_q\}
       \cup\{sI:s\in\F_q\setminus\{0,1,-1\}\}
\label{For_D}
\end{equation}
is a dominating set of $\tau(\Mat)$ and has cardinality $2q-2$.
Moreover, the matrices $A\in\Mat$ for which $A+B$ is invertible
for every $B\in\D$ are exactly
\[
 \mathcal U_q=\{I+tE_{12}:t\in\F_q\}.
\]
\end{proposition}
\begin{proof}
Notice that the characteristic is odd, so the elements $0$, $1$, and $-1$ are distinct. The three sets in (\ref{For_D}) are pairwise disjoint, so \( \left\vert \D\right\vert =1+q+(q-3)=2q-2
\).
Let $A=\begin{pmatrix}a&b\\c&d\end{pmatrix}\in\Mat\setminus\D$. If $A$ is singular, then $A\ne0$ and $A$ is adjacent to $0\in\D$. We may therefore assume that $A$ is invertible.

Assume first that $c\ne0$. For every $t\in\F_q$,
\[
 \det(A+I+tE_{12})=(a+1)(d+1)-c(b+t).
\]
Taking \( t=c^{-1}\big((a+1)(d+1)-cb\big)
\)
makes this determinant zero. Since $I+tE_{12}\in\D$ and $A\not \in \D$, these are distinct adjacent vertices.

Now assume that $c=0$. Invertibility gives $a,d\ne0$. Notice that $\D$ contains $sI$ for every $s\in\F_q\setminus\{-1\}$. If $a\ne1$, choose $s=-a$. Then $s\not =-1$ and
\[
 \det(A+sI)=(a+s)(d+s)=0.
\]
Similarly, if $d\ne1$, choose $s=-d$. If neither choice applies, then $a=d=1$ and $A=I+bE_{12}\in\D$, contrary to the choice of $A$. Every vertex outside $\D$ is therefore adjacent to a member of $\D$.

A vertex in $\D$ need not have a neighbour in $\D$. We now identify
the matrices whose sums with every member of $\D$ are invertible.
If $A$ is such a matrix, then $A$ is invertible because $0\in\D$.
The determinant calculation above forces $c=0$, and the scalar
matrices in $\D$ then force $a=d=1$, so $A\in\mathcal U_q$.
Conversely, let $A=I+tE_{12}$. Its sum with $0$ is invertible.
Its sum with $I+uE_{12}$ has determinant $4\ne0$. Its sum with any
other scalar matrix $sI\in\D$ has determinant $(1+s)^2\ne0$, since
$s\not=-1$. Thus all the sums are invertible.
\end{proof}

\begin{proof}[Proof of Theorem~\ref{Tm_main}(a)]
The upper bound $\gamma(\tau(\Mat))\le2q-2$ is
Proposition~\ref{Pr_upper}. We prove the lower bound
$\gamma(\tau(\Mat))\ge q+1$ for every prime power $q\ge3$, odd or even.
It suffices to show that no set $D$ of $q$ matrices dominates.
After interchanging the two rows of every matrix if necessary,
we may assume that two members of $D$ have distinct first rows
$a_1,a_2$. Let $A$ be the set of first rows of the members of $D$,
and let $\ell$ be the affine line through $-a_1$ and $-a_2$.
Since
\[
|\ell\cup A\cup(-A)|\le q+q+(q-2)=3q-2<q^2,
\]
we may choose $r\in\F_q^2$ outside this union.

Consider the matrices $X$ with first row $r$ and second row $s$.
For a member of $D$ with rows $a,b$, the sum with $X$ is singular
precisely when
\[
s\in -b+\langle r+a\rangle.
\]
This is an affine line since $r\notin -A$. Since $r\notin\ell$, the points $-a_1,-a_2,r$ are not collinear.
Hence the direction vectors $r+a_1$ and $r+a_2$ are linearly
independent, so the two corresponding lines are nonparallel.
Each of the $q$ members of $D$ thus rules out one affine line of
$q$ possible second rows. Since two of these lines are nonparallel,
they intersect, so their union together with the remaining lines
contains at most $q^2-1$ points.
Choose $s$ outside their union. Then $X$ has invertible sum with
every member of $D$, and $X\notin D$ since $r\notin A$.
Hence $D$ does not dominate.

Finally, since $2q-2<2q-1$, the upper bound contradicts the value
predicted by~\eqref{For_conjecture} for $R=M_2(\F_q)$.
\end{proof}

In the smallest case $q=3$, both bounds equal $4$, so
$\gamma(\tau(M_2(\F_3)))=4$, and
$\D=\{0,I,I+E_{12},I-E_{12}\}$ is a minimum dominating set.

\section{Total domination and singular-sum coverings}\label{Sec_covering}

For a dominating set to be total, each of its vertices must also have a neighbour in the set. A \emph{total dominating set}
of a graph $G$ is a set $T$ such that every vertex of $G$ has a distinct
neighbour in $T$. Its minimum cardinality is denoted by $\gamma_t(G)$,
provided such a set exists.

To study this parameter, we introduce an intermediate covering number.
For a finite nonzero ring $R$ with identity, define
\begin{equation}\label{For_kappa}
\kappa(R)=\min\bigl\{|D|:D\subseteq R,\
\text{for every }x\in R\text{ there is }d\in D
\text{ with }x+d\in Z(R)\bigr\}.
\end{equation}
We call such a set $D$ a \emph{zero-divisor-sum covering set}, or a
\emph{singular-sum covering set} in a matrix ring over a field.
The condition applies also to selected elements, but permits $d=x$
when $2x\in Z(R)$. Equivalently, the translates $-d+Z(R)$, for $d\in D$,
cover $R$. We shall prove that this covering number equals the total
domination number for every finite nonlocal ring of odd
order.

\begin{proposition}
\label{Prop_three}
\begin{enumerate}
\item[(a)] For a finite nonzero ring $R$,
\[
 \gamma(\tau(R))\leq\kappa(R)\leq\gamma_t(\tau(R)),
\]
where the second inequality holds whenever $\tau(R)$ has no isolated vertex. If $R$ has characteristic two, then $\gamma(\tau(R))=\kappa(R)$.
\item[(b)] For every prime power $q$ and $n\geq2$,
\[
 \gamma(\tau(M_n(\F_q)))\leq\kappa(M_n(\F_q))
 \leq\gamma_t(\tau(M_n(\F_q)))\leq n(q-1)+1.
\]
\end{enumerate}
\end{proposition}
\begin{proof}
(a) A set satisfying (\ref{For_kappa}) provides a distinct neighbour for every vertex outside that set, so it dominates $\tau(R)$. Every total dominating set satisfies (\ref{For_kappa}), giving the second inequality. In characteristic two, a dominating set automatically satisfies the covering condition at each selected element $d$, because $d+d=0\in Z(R)$.

(b) Consider the canonical set from \cite[Lemma~4.2]{DO},
\[
 D_0=\{0\}\cup\{xE_{1j}:x\in\F_q^\times,\ 1\leq j\leq n\}.
\]
If $A$ is singular, $A+0$ is singular. If $A$ is invertible, some cofactor $C_{1j}$ is nonzero, and
$\det(A+xE_{1j})=\det A+xC_{1j}$ vanishes for an appropriate nonzero $x$. Thus $D_0$ satisfies the covering condition. Moreover, every sum of two members of $D_0$ has rank at most one, hence is singular. Since $\left\vert D_0\right\vert \geq3$, the selected vertices have distinct neighbours in $D_0$. It follows that $D_0$ is a total dominating set. Now apply part~(a).
\end{proof}

The parameters in (a) need not coincide in general. For example, $\kappa(\F_3)=3$, whereas $\tau(\F_3)$ has an isolated vertex and hence has no total dominating set.

For nonlocal rings of odd order, covering and total domination
have the same minimum cardinality.

\begin{theorem}\label{Tm_shift}
Let $R$ be a finite nonlocal ring of odd order. Every
zero-divisor-sum covering set in $R$ has an additive translate
that is a total dominating set of $\tau(R)$. Consequently,
\[
\gamma_t(\tau(R))=\kappa(R).
\]
\end{theorem}

\begin{proof}
We first show that every element of $R$ is a sum of two
zero-divisors. Write
\[
R/J(R)\cong\prod_{i=1}^k M_{n_i}(\F_{q_i}).
\]
If $k\geq2$, every element of this product is a sum of two
nonunits: split its coordinates into two nonempty groups
and set the remaining coordinates of each summand equal
to zero. If $k=1$, nonlocality implies $n_1\geq2$.
Every matrix in $M_{n_1}(\F_{q_1})$ is a sum of two singular
matrices, obtained by separating its last column from the
other columns.

For any $a\in R$, lift such a decomposition of $a+J(R)$
by choosing a lift $u$ of the first summand and setting
$v=a-u$. Both $u$ and $v$ are nonunits, because their images
modulo $J(R)$ are nonunits. Since $R$ is finite, nonunits
are precisely zero-divisors. Thus
\[
R=Z(R)+Z(R).
\]

Let $D$ be a covering set and put $N=|R|$. A single
translate of the proper subset $Z(R)$ cannot cover $R$,
so $|D|\geq2$. For each $d\in D$, define
\[
P_d=\{x\in R:x+d\in Z(R),\
x+e\notin Z(R)\text{ for every }e\in D\setminus\{d\}\}.
\]
These private target sets are pairwise disjoint, and
\begin{equation}\label{For_private}
\sum_{d\in D}|P_d|<N.
\end{equation}
Indeed, choose distinct $d,e\in D$ and write
$e-d=u+v$ with $u,v\in Z(R)$. Then $x=-u-d$ satisfies
\[
x+d=-u\in Z(R),\qquad x+e=v\in Z(R).
\]
Hence $x$ belongs to none of the private target sets.

For every $h\in R$, the translate $D+h$ is again a covering
set: apply the covering property of $D$ to $x+h$ for any
$x\in R$. A selected vertex $d+h$ has no distinct neighbour
in $D+h$ precisely when \(d+2h\in P_d
\).
Indeed, its sums with $e+h$, for $e\ne d$, are
$(d+2h)+e$. If all these are units, the covering property
forces $(d+2h)+d$ to be a zero-divisor.

Since the additive group of $R$ has odd order, the map
$h\mapsto2h$ is bijective. Thus exactly $|P_d|$ shifts are
bad for a given $d$. By~\eqref{For_private}, fewer than $N$
shifts are bad for at least one member of $D$. Choose a
shift that is bad for none. Every selected vertex then
has a distinct selected neighbour, and every vertex
outside $D+h$ has one by the covering property.

Applying this to a minimum covering set gives
$\gamma_t(\tau(R))\leq\kappa(R)$.
The reverse inequality follows from
Proposition~\ref{Prop_three}.
\end{proof}

The covering condition has the structural properties needed to reduce finite rings to full matrix rings. The formulas for $\kappa$ in the following proposition also follow
from Handelman's Lemmas~A.1 and~A.2, stated in terms of
property $((k))$~\cite{Handelman}.

\begin{proposition}\label{Pr_structure}
Let $R,S$ be finite nonzero rings with identity, and let
$J=J(R)$. Then
\[
\kappa(R)=\kappa(R/J)
\]
and
\(
\gamma_t(\tau(R\times S))
=\kappa(R\times S)
=\min\{\kappa(R),\kappa(S)\}
\).
\end{proposition}
\begin{proof}
In a finite ring, nonunits are exactly zero-divisors, and an element is a unit if and only if its image modulo the Jacobson radical is a unit. Therefore
\[
 x+d\in Z(R)\quad\Longleftrightarrow\quad
 (x+J)+(d+J)\in Z(R/J).
\]
A covering set projects to a covering set of no greater cardinality. Conversely, choosing one representative of each element of a covering set in $R/J$ gives a covering set in $R$, because the displayed condition is valid for every $x$. This proves the first identity.

In $R\times S$, an element is a zero-divisor precisely when
at least one of its coordinates is a zero-divisor. Suppose,
after interchanging the factors if necessary, that
$\kappa(R)\leq\kappa(S)$, and let $D$ be a minimum covering
set in $R$. Then $D\times\{0\}$ covers $R\times S$.
Moreover, $|D|\geq2$, since a single translate of the proper
subset $Z(R)$ cannot cover $R$. Any two distinct members of
$D\times\{0\}$ are adjacent, since the second coordinate of
their sum is zero. Thus $D\times\{0\}$ is a total dominating
set, and Proposition~\ref{Prop_three} gives
\[
\kappa(R\times S)\leq\gamma_t(\tau(R\times S))
\leq\min\{\kappa(R),\kappa(S)\}.
\]

For the reverse bound, let $\{(a_i,b_i):1\leq i\leq m\}$ cover $R\times S$. If $\{a_i\}$ covers $R$, then $\kappa(R)\leq m$. Otherwise choose $x\in R$ for which every $x+a_i$ is a unit. For each $y\in S$, covering $(x,y)$ forces $y+b_i\in Z(S)$ for some $i$. Thus $\{b_i\}$ covers $S$ and $\kappa(S)\leq m$. In either case, $\min\{\kappa(R),\kappa(S)\}\leq m$.
\end{proof}

A nonzero radical also determines the relation between
 domination and covering.

\begin{theorem}\label{Tm_radical_domination}
Let $R$ be a finite nonzero ring with $J=J(R)\ne0$. Then
\[
\gamma(\tau(R))=\kappa(R)=\kappa(R/J).
\]
If, in addition, $R$ is nonlocal of odd order, then also
$\gamma_t(\tau(R))=\kappa(R)$, so all four quantities coincide.
\end{theorem}

\begin{proof}
Let $\pi:R\to R/J$ be the quotient map, let $D$ be a
dominating set of $\tau(R)$, and put $E=\pi(D)$. Let $U$
be the set of elements of $R/J$ uncovered by $E$, that is,
\[
U=\{u\in R/J: u+e\text{ is a unit for every }e\in E\}.
\]
If $u\in U$, then every element of $\pi^{-1}(u)$ has unit
sum with every member of $D$. Since $D$ dominates $\tau(R)$,
the whole fibre $\pi^{-1}(u)$ must therefore lie in $D$.
In particular, $U\subseteq E$.

Each fibre over $E\setminus U$ contributes at least one
element to $D$, and each fibre over $U$ contributes all
$|J|$ elements. Since $|J|\geq2$, it follows that
\[
|D|\geq |E\setminus U|+|J|\,|U|
     = |E|+(|J|-1)|U|
     \geq |E|+|U|.
\]
Now $E\cup(-U)$ covers $R/J$. Indeed, $E$ already covers
every element outside $U$, and each $u\in U$ is covered
by $-u$. Consequently,
\[
\kappa(R/J)\leq |E\cup(-U)|\leq |E|+|U|\leq |D|.
\]
Taking the minimum over all dominating sets $D$, we obtain
$\kappa(R/J)\leq\gamma(\tau(R))$. The reverse inequalities
follow from Proposition~\ref{Prop_three} and
Proposition~\ref{Pr_structure}, which give
\(\gamma(\tau(R))\leq\kappa(R)=\kappa(R/J)\).

If $R$ is moreover nonlocal of odd order, then
Theorem~\ref{Tm_shift} gives $\gamma_t(\tau(R))=\kappa(R)$.
\end{proof}

The following example shows that the three domination and covering
parameters can differ when \(R\) is local or when \(J(R)=0\).

\begin{example}\label{Ex_hypotheses}
\emph{Local rings.} Let $R$ be a finite local ring of odd order with
$J=J(R)\ne0$, and put $q=|R/J|$. Since $R/J\cong\F_q$,
Theorem~\ref{Tm_radical_domination} gives
$\gamma(\tau(R))=\kappa(R)=q$. On the other hand, put $a=|J|>1$.
Two distinct elements of $R$ are adjacent
precisely when their images in $R/J$ sum to zero.
The coset $J$ therefore induces a complete graph $K_a$,
each pair of nonzero opposite cosets induces a complete
bipartite graph $K_{a,a}$, and there are no edges between
these components. Thus
\[
\tau(R)\cong
K_a\ \sqcup\
\bigsqcup_{i=1}^{(q-1)/2}K_{a,a}.
\]
Each component has total domination number $2$, so
\(
\gamma_t(\tau(R))
=2+2\frac{q-1}{2}=q+1>\kappa(R)
\).
Hence the nonlocal hypothesis in Theorem~\ref{Tm_shift} is needed.

\emph{Semisimple rings.} Domination and total domination differ
for \(M_2(\mathbb F_3)\):
\[
 \gamma(\tau(M_2(\F_3)))=4,\qquad
 \kappa(M_2(\F_3))=\gamma_t(\tau(M_2(\F_3)))=5.
\]
The equality $\gamma(\tau(M_2(\F_3)))=4$ follows from
Theorem~\ref{Tm_main}(a). The geometric bound in
Corollary~\ref{cor:sqrt} below gives
$\kappa(M_2(\F_3))\geq5$.
The canonical set $\{0,E_{11},-E_{11},E_{12},-E_{12}\}$ is a total
dominating set by Proposition~\ref{Prop_three}(b), proving the two
remaining equalities. The value $\kappa(M_2(\F_3))=5$ is also noted
in~\cite[immediately after Lemma~A.16]{Handelman}.
In contrast, let $A=\F_3[\varepsilon]/(\varepsilon^2)$
and $R=M_2(A)$. Then
\[
J(R)=M_2((\varepsilon))\ne0,
\qquad R/J(R)\cong M_2(\F_3).
\]
Theorem~\ref{Tm_radical_domination} therefore gives
\[
\gamma(\tau(R))=\kappa(R)=\gamma_t(\tau(R))=5.
\]
Thus the domination number is $4$ for the total graph of the semisimple quotient and $5$ for that of $R$.
\end{example}

\section{Geometric lower bounds for coverings}\label{Sec_geometric}
For $q\ge3$, let $b(q)$ be the smallest size of a point set in
$\PG(2,q)$ that meets every line and contains no line. Such a set is
called a \emph{nontrivial blocking set}.

View $\PG(3,q)$ as the projective space associated with
$\F_q^2\oplus\F_q^2$. Let $L=L_\infty$ be the line
corresponding to $\{(0,w):w\in\F_q^2\}$. For $X\in\Mat$,
let $L_X$ be the line corresponding to the subspace
\[
\{(v,Xv):v\in\F_q^2\}.
\]
Every line skew to $L$ is $L_X$ for a unique $X\in\Mat$,
since projection onto the first factor restricts to an
isomorphism on its corresponding two-dimensional subspace.
Moreover, $L_X\cap L_A\ne\varnothing$ if and only if
$\det(X-A)=0$, since a common point corresponds to a nonzero
vector $v$ satisfying $Xv=Av$. This includes the case $X=A$.
Thus $D\subseteq\Mat$ is a singular-sum covering set if and only
if every line skew to $L$ has a point in common with some member
of $C=\{L_{-d}:d\in D\}$.

For the remainder of this section, assume that $D$ is such a
covering set. Write $m=|C|$ and $S=\bigcup_{c\in C}c$.

\begin{lemma}\label{lem:missingpoint}
If $m\le2q-2$ and $g\ne L$ meets $L$ at $R$, then $g\setminus\{R\}$
contains a point outside $S$.
\end{lemma}
\begin{proof}
Suppose that all $q$ points of $g\setminus\{R\}$ belong to $S$, and
put $\rho=\langle g,L\rangle$. Each member of $C$ meets $\rho$ in exactly one point, since a line
contained in $\rho$ would meet $L$.
At least $q$ members therefore meet $g$, leaving at most $q-2$ points
of $S\cap\rho$ outside $g$. Of the $q-1$ lines of $\rho$ through $R$
other than $g$ and $L$, at least one, say $g'$, misses $S$.
Fix $Z\in g'\setminus\{R\}$. Every line through $Z$ other than $g'$
meets $S$. Indeed, a line through $Z$ in $\rho$ other than $g'$ meets
$g$ away from $R$, and a line through $Z$ outside $\rho$ is skew to $L$.

Project from $Z$ onto a plane $\pi$ avoiding $Z$. For each $c\in C$
the plane $\langle Z,c\rangle$ cuts $\pi$ in a line. A point $X\in\pi$ lies on one of the projected lines if and only if
$ZX$ meets $S$. For $P=g'\cap\pi$, we have $ZP=g'$, which is
disjoint from $S$, whereas every other line through $Z$ meets $S$.
Thus the projected lines cover precisely $\pi\setminus\{P\}$.
In the dual plane, with $P^*$ as the line at infinity, they
correspond to at most $m$ affine points meeting every affine line.
The affine blocking-set theorem~\cite[p.~251, Theorem]{BS} gives
$m\ge2q-1$, a contradiction.
\end{proof}

\begin{lemma}\label{L:general}
For every prime power $q\ge3$,
\[
\kappa(\Mat)\ge\min\{2q-1,b(q)\}.
\]
\end{lemma}
\begin{proof}
Suppose that $m\le2q-2$. If the members of $C$ are pairwise skew,
then $C\cup\{L\}$ is a maximal partial line spread, since every
line meets one of its members. Its size is
$m+1\le2q-1<q^2+1$, so it is not a spread. Glynn's lower
bound~\cite{Glynn}, recalled in
\cite[p.~48, Theorem~1.2]{DeBoeck}, gives $m+1\ge2q$,
a contradiction.

Choose two members meeting at $Y$, and a plane $\rho$ through $L$
avoiding $Y$. The set $S\cap\rho$ has at most $m$ points, all outside
$L$. By \cite[p.~251, Theorem]{BS} applied
with $k=2$, a point set meeting every affine line of $\rho\setminus L$
has at least $2q-1$ points. Since $|S\cap\rho|\le m\le 2q-2$,
some affine line misses $S\cap\rho$. Its projective completion $n$
is disjoint from $S$ because $S\cap L=\varnothing$. Put $\sigma=\langle n,Y\rangle$
and $Q=n\cap L$. Then $\sigma\cap L=\{Q\}$, and $\sigma$ contains
no member of $C$, since such a member would meet $n$.

Put $B=(S\cap\sigma)\cup\{Q\}$. Since $n$ misses $S$ and
$Q\in n$, we have $B\cap n=\{Q\}$. We claim that $B$ meets
every line of $\sigma$.
Indeed, lines through $Q$ meet $B$ there, and every other line is
skew to $L$ and hence meets $S$. Moreover, $B$ contains no line:
a line of $\sigma$ avoiding $Q$ meets $n$ outside $B$, while a line
through $Q$ has a point outside $S\cup\{Q\}$ by
Lemma~\ref{lem:missingpoint}. Thus $|B|\ge b(q)$. Each member of
$C$ contributes one point to $S\cap\sigma$, and the two chosen
members contribute the same point $Y$. Hence $|B|\le m$, as required.
\end{proof} 

For $q\ge7$, a closer look at the case of equality yields a strict
lower bound in terms of $b(q)$.

\begin{theorem}\label{thm:strict-blocking}
For every prime power $q\ge7$,
\[
\kappa(\Mat)\ge b(q)+1.
\]
\end{theorem}
\begin{proof}
We first record that $b(q)\le3(q+1)/2$. For a function
$f:\F_q\to\F_q$, take its affine graph together with the points
at infinity of its secant directions. This set meets every line:
in a missing direction the $q$ graph points lie on distinct parallel
lines. It contains no line if $f$ is nonaffine, since the vertical
direction is missing and no affine line contains all $q$ graph points.
For odd $q$, take $f(0)=0$, with $f(x)=x$ on nonzero squares and
$f(x)=-x$ on nonsquares. Besides $1$ and $-1$, the slopes have the
form $(1+t)/(1-t)$ with $t$ a nonsquare, so there are at most
$(q+3)/2$ directions. For even $q$, take the absolute trace
$f(x)=\operatorname{Tr}_{\F_q/\F_2}(x)$. Its nonzero slopes are
$1/h$ with $\operatorname{Tr}(h)=1$, so there are at most $q/2+1$
directions including slope zero. Both functions are nonaffine for
$q\ge7$, proving the claimed bound.

Since $b(q)\le3(q+1)/2\le2q-2$, Lemma~\ref{L:general}
reduces the proof to excluding a cover with $m=b(q)$ lines.
In the notation above, call a line $n\ne L$ free if it meets $L$
and misses $S$. Every plane $\rho$ through $L$ contains a free
line, by the affine blocking-set bound applied to $S\cap\rho$.
For any plane $\sigma$ through such an $n$, other than
$\langle n,L\rangle$, the proof of Lemma~\ref{L:general}
shows that $(S\cap\sigma)\cup\{n\cap L\}$ is a nontrivial
blocking set. Since $n$ is disjoint from $S$, no member of $C$ lies in $\sigma$.
The $m$ members therefore meet $\sigma$ in $m$ points counted with
multiplicity, of which at least $b(q)-1=m-1$ are distinct. Thus at
most one point occurs twice, and none occurs three times.

For any $Y\in S$, choose a plane through $L$ avoiding $Y$ and
a free line $n$ in it. Applying the preceding observation to
$\langle n,Y\rangle$ shows that at most two members pass through
$Y$. Under projective duality, the dual family consists of $m=b(q)$
lines skew to the dual of $L$, and every line skew to the dual
of $L$ meets at least one member of this family. The dual family therefore satisfies the same hypotheses relative to
the dual of $L$. By the point bound just proved, at most two dual
members pass through any point. Since the members of $C$ contained
in a plane correspond to the dual members passing through the dual
point of that plane, at most two members of $C$ lie in any plane.
Hence no three members are pairwise intersecting, since three
pairwise intersecting lines are concurrent or coplanar.

Let $J$ count the intersecting pairs in $C$. Their intersection
points are distinct. For each plane $\rho$ through $L$, let
$J_\rho$ count those points in $\rho$. Choose a free line in $\rho$. Each of the $q$ other planes
$\sigma$ through this line contains at most one such point:
two distinct intersection points would give
$|S\cap\sigma|\le m-2$, whereas $|S\cap\sigma|\ge m-1$.
Thus $J-J_\rho\le q$. Each intersection point belongs to exactly
one plane through $L$, so summing over all $q+1$ planes gives
$J\le q+1$.

Put $M=C\cup\{L\}$ and $v=|M|=m+1$. For $u\notin M$ let
$i(u)$ count the members of $M$ meeting $u$. Every such line has
$i(u)\ge1$. Write $(a)_k=a(a-1)\cdots(a-k+1)$, with $(a)_0=1$,
and put $S_k=\sum_{u\notin M}(i(u))_k$.  Let $t_j$ count the unordered triples in $M$ with $j$
intersecting pairs. Since no three members of $C$ are pairwise
intersecting and $L$ is skew to every member of $C$, we have
$t_3=0$. Put
$P=\sum_{w\in M}(d_w)_3$, where $d_w$ counts the other members
meeting $w$. Counting ordered incidences gives
\begin{align*}
S_0&=(q^2+1)(q^2+q+1)-v,\\
S_1&=vq(q+1)^2-2J,\\
S_2&=(v)_2(q+1)^2+2(q^2-q-2)J-2t_2,\\
S_3&=(v)_3(q+1)+6q(v-2)J-6(q+1)t_2-P.
\end{align*}
Indeed, a line meets $q(q+1)^2$ other lines. Two skew lines
have $(q+1)^2$ common transversals. For two intersecting lines,
the transversals lie in their plane or pass through their common
point, giving $2q^2+q-1$ after the two lines themselves are
excluded. A member of $M$ counted as a transversal to two other
members contributes to the correction $2t_2$.

For triples with zero, one, or two intersecting pairs, the respective
numbers of transversals outside the triple are $q+1,2q+1,2q$.
For three skew lines, each point of the first lies on a unique
transversal: the plane it spans with the second meets the third
in one point. If just two lines intersect, their plane meets the
third line at a point $T$. The transversals form the pencil in
that plane through $T$ and the pencil joining their intersection
point to the third line, with one line in common. If two pairs
intersect, the common member meets the two skew members at
points $A,B$. A transversal other than the common member joins
$A$ to a point other than $B$ on the opposite line, or $B$ to a
point other than $A$ on the other line. Now use
$t_0+t_1+t_2=\binom v3$ and $t_1+2t_2=(v-2)J$, and subtract
$P$ to exclude transversals in $M$. This gives the formula for $S_3$.

Since $(i-1)(i-3)(i-4)\ge0$ for every integer $i\ge1$, we have
\begin{align*}
0&\le S_3-5S_2+12S_1-12S_0\\
 &= F_q(v)+\bigl(6q(v-2)-10(q^2-q-2)-24\bigr)J
       -(6q-4)t_2-P,
\end{align*}
where
\[
F_q(v)=(v)_3(q+1)-5(v)_2(q+1)^2
 +12vq(q+1)^2-12\bigl((q^2+1)(q^2+q+1)-v\bigr).
\]
Put $V=(3q+5)/2$, so $v\le V$. On $[0,V]$,
$F_q''(v)\le F_q''(V)=-(q+1)^2$ and
$F_q'(V)=(5q+9)(3q^2-10q-1)/4>0$. Hence $F_q(v)\le F_q(V)$.
The coefficient of $J$ is at most $-q^2+13q-4\le38$.
Using $J\le q+1$ and writing $x=q-7\ge0$, we obtain
\[
0\le 8\bigl(F_q(V)+38(q+1)\bigr)
=-352-4736x-2028x^2-300x^3-15x^4<0.
\]
This excludes $m=b(q)$. 
\end{proof}

Combining Lemma~\ref{L:general} and Theorem~\ref{thm:strict-blocking}
with known lower bounds for $b(q)$ gives explicit estimates for
$\kappa(\Mat)$.

\begin{corollary}\label{cor:sqrt}
For every prime power $q\ge3$,
\[
\kappa(\Mat)\ge
\begin{cases}
2q-1, & q\in\{3,4,5\},\\
q+\lceil\sqrt q\rceil+2, & q\ge7.
\end{cases}
\]
For every prime $p\ge7$,
\[
\kappa(M_2(\F_p))\ge(3p+5)/2.
\]
\end{corollary}
\begin{proof}
Bruen's bound~\cite{Bruen}, recalled in~\cite[p.~111]{Blokhuis}, gives
$b(q)\ge q+\lceil\sqrt q\rceil+1$. For odd primes $p$,
Blokhuis's theorem~\cite[p.~111, Theorem]{Blokhuis} gives
$b(p)\ge 3(p+1)/2$.
For $q=3,4,5$, the first bound and Lemma~\ref{L:general}
give $\kappa(\Mat)\ge2q-1$. For $q\ge7$, apply
Theorem~\ref{thm:strict-blocking}.
\end{proof}

\section{Domination and uncovered matrices}\label{Sec_fibres}

For $D\subseteq\Mat$, write
\[
 W(D)=\{X\in\Mat:\det(X+d)\ne0\text{ for all }d\in D\}.
\]
For example, for odd $q$, Proposition~\ref{Pr_upper} gives
$W(\D)=\{I+tE_{12}:t\in\F_q\}$. If $D$ dominates, then $W(D)\subseteq D$. We show that a dominating set
with $W(D)\ne\varnothing$ must have at least $2q-2$ elements for odd $q\ge3$.

For a nonzero row vector $\varphi\in\F_q^2$ and a row
$r\in\F_q^2$, the fibre $\{x\in\Mat:\varphi x=r\}$
is an affine plane.

\begin{lemma}\label{L:fibre}
Let $D\subseteq\Mat$ be any set with $m$ elements, $w_0\in W(D)$, $\varphi\ne0$ and $\Pi=\{x:\varphi x=\varphi w_0\}$. Then
\[
|W(D)\cap\Pi|\ \ge\ 2q-1-m .
\]
\end{lemma}
\begin{proof}
Choose $P\in\GL_2(\F_q)$ with first row
$\varphi$. Since
\[
\det(Px+Pd)=\det(P)\det(x+d),
\]
the map $x\mapsto Px$ sends $W(D)$ bijectively onto $W(PD)$
and sends $\Pi$ onto the set of matrices with first row
$r=\varphi w_0$. Thus, after this change of coordinates, we
may assume that $\varphi=e_1^T$.

Write $r=(r_1,r_2)$. Each matrix in $\Pi$ has the form
\[
x(y)=
\begin{pmatrix}
r_1 & r_2\\
y_1 & y_2
\end{pmatrix},
\qquad y=(y_1,y_2)\in\F_q^2.
\]
For each $d=(d_{ij})\in D$, the polynomial
\[
L_d(y)=\det(x(y)+d)
      =(r_1+d_{11})(y_2+d_{22})
       -(r_2+d_{12})(y_1+d_{21})
\]
has degree at most one in $y_1,y_2$. Hence
\[
G(y)=\prod_{d\in D}L_d(y)
\]
has total degree at most $m$. By definition, $G(y)\ne0$ if
and only if $x(y)\in W(D)\cap\Pi$. In particular, $G$ is
nonzero at the point corresponding to $w_0$.

By the Alon--F\"uredi theorem~\cite[p.~81, Theorem~5]{AF},
the number of points of $\F_q^2$ at which $G$ is nonzero
is at least the minimum of $ab$ over all integers
$1\leq a,b\leq q$ satisfying $a+b\geq 2q-\deg G$.
For every such pair,
\[
ab=(a-1)(b-1)+a+b-1
   \geq a+b-1
   \geq 2q-1-\deg G
   \geq 2q-1-m.
\]
Since $y\mapsto x(y)$ is a bijection from $\F_q^2$ onto
$\Pi$, this gives $|W(D)\cap\Pi|\geq 2q-1-m$.
\end{proof}

\begin{proof}[Proof of Theorem~\ref{Tm_main}(b)]
Proposition~\ref{Prop_three} and Proposition~\ref{Pr_upper} give the
upper bound. Let $D$ be a minimum dominating set and write $W=W(D)$.
If $W=\varnothing$, then $D$ is a covering set and
$|D|\ge\kappa(\Mat)$. It remains to show that
$|D|\ge2q-2$ when $W\ne\varnothing$.

Suppose, to the contrary, that $W\ne\varnothing$ and
$m=|D|\le2q-3$. Since $D$
dominates, $W\subseteq D$. Join two elements of $W$ when their
difference has rank one.

For each $u\in W$ and each left kernel direction, the corresponding
fibre through $u$ contains at least $2q-1-m\ge2$ elements of $W$,
by Lemma~\ref{L:fibre}. Thus $u$ has a neighbour in every left
kernel direction and, by transposition, in every right kernel
direction. In particular, every vertex has at least $q+1$ neighbours.

Fix $w_0\in W$, let $N_0$ be its neighbourhood, and put $k=|N_0|$.
Label each $w\in N_0$ by the left and right kernels of $w-w_0$.
All $q+1$ left labels occur. Since
\[
k\le |W|-1\le 2q-4<2(q+1),
\]
some left label occurs only once. Let $w$ be the vertex with
this label.

Two distinct vertices of $N_0$ are adjacent precisely when they
share a left or a right label. Indeed, the difference of two
rank-one matrices $ab^T$ and $cd^T$ is singular precisely when
$a,c$ are dependent or $b,d$ are dependent. Hence every neighbour
of $w$ in $N_0$ shares its right label. The other $q$ right labels
all occur, so at most $k-q-1$ vertices of $N_0$ are neighbours
of $w$. There are $|W|-k$ vertices outside $N_0$. Consequently,
\[
\deg(w)\le(k-q-1)+(|W|-k)=|W|-q-1\le q-4,
\]
contrary to $\deg(w)\ge q+1$.
\end{proof}

\begin{proof}[Proof of Theorem~\ref{Tm_main}(c)]
Corollary~\ref{cor:sqrt} gives $\kappa(\Mat)\ge 2q-1$ for
$q=3,4,5$. For $q=7$, its prime-field bound gives
$\kappa(M_2(\F_7))\ge 13$. Thus, for all four field orders,
Proposition~\ref{Prop_three}(b) yields
\[
2q-1\le\kappa(\Mat)\le\gamma_t(\tau(\Mat))\le2q-1.
\]
For odd $q$, part~(b) then gives $\gamma=2q-2$, and for $q=4$
Proposition~\ref{Prop_three}(a) gives $\gamma=\kappa=7$. 
\end{proof}

\section{Matrices of size at least three}\label{Sec_n3}

For $n\ge2$ and $D\subseteq M_n(\F_q)$, extend the notation of
Section~\ref{Sec_fibres} by setting
\[
W(D)=\{X\in M_n(\F_q):\det(X+d)\ne0\text{ for all }d\in D\}.
\]
As before, $D$ is a covering set exactly when $W(D)=\varnothing$,
and it is a dominating set exactly when $W(D)\subseteq D$.
The next two lemmas hold for $n\ge2$.

 If $w_0\in W(D)$ and
$d\in D$, then $a_d=w_0+d$ is invertible, and for $u,r\in\F_q^n$
\begin{equation}\label{For_rank_one}
 \det(w_0+ur^T+d)=\det(a_d)\bigl(1+r^Ta_d^{-1}u\bigr).
\end{equation}
We use the following covering theorem of
Jamison~\cite[p.~257, Theorem~$1'$]{Jamison}.
\emph{If a family of affine
hyperplanes of $\F_q^n$, none of which contains $0$, covers
$\F_q^n\setminus\{0\}$, then the family has at least $n(q-1)$ members.}

\begin{lemma}\label{L:coincide}
Let $n\ge2$, let $D\subseteq M_n(\F_q)$ with $|D|\le n(q-1)$, let
$w_0\in W(D)$, and let $d\ne d'$ be members of $D$.
\begin{enumerate}
\item[(a)] If $v\ne0$ and $(d-d')v=0$, put $u=(w_0+d)v$. Then $u\ne0$,
and $w_0+ur^T\in W(D)$ for some $r\ne0$.
\item[(b)] If $y\ne0$ and $y^T(d-d')=0$, put $r^T=y^T(w_0+d)$. Then
$r\ne0$, and $w_0+ur^T\in W(D)$ for some $u\ne0$.
\end{enumerate}
\end{lemma}

\begin{proof}
(a) Since $(d-d')v=0$, we have $u=a_dv=a_{d'}v$, and $u\ne0$ because
$a_d$ is invertible. For $e\in D$ put $c_e=a_e^{-1}u\ne0$. By
\eqref{For_rank_one}, the matrix $w_0+ur^T+e$ is singular precisely when
$r$ lies on the affine hyperplane $H_e=\{r:r^Tc_e=-1\}$, which does not
contain $0$. Since $c_d=v=c_{d'}$, we have $H_d=H_{d'}$, so among the
hyperplanes $H_e$, $e\in D$, there are at most
$|D|-1<n(q-1)$ distinct ones. By the covering theorem above, some $r\ne0$
lies on none of them, and then $w_0+ur^T\in W(D)$.

(b) Since $W(D^T)=W(D)^T$, this follows by applying~(a) to $w_0^T$
and the set $D^T$.
\end{proof}

For a subspace $R$ of the row space $\F_q^n$, let $L_R$ be the set of
matrices all of whose rows lie in $R$.

\begin{lemma}\label{L:sum}
Let $n\ge2$, let $D\subseteq M_n(\F_q)$ with $|D|<(n+1)(q-1)$, let
$X_0\in M_n(\F_q)$, and let $R$ be a $2$-dimensional subspace of
$\F_q^n$. Then
\[
 \sum_{Y\in L_R}\ \prod_{d\in D}\det(X_0+Y+d)=0 .
\]
Consequently, $W(D)\cap(X_0+L_R)$ does not consist of exactly one
matrix.
\end{lemma}

\begin{proof}
Choose $Q\in\GL_n(\F_q)$ such that $RQ$ is spanned by the first two
standard row vectors. Replacing $X_0$ and $D$ by $X_0Q$ and $DQ$
multiplies every determinant by $\det Q$ and maps $L_R$ onto $L_{RQ}$.
We may therefore assume that $L_R$ consists of the matrices
$Y=[\,y_1\ y_2\ 0\ \cdots\ 0\,]$ with columns $y_1,y_2\in\F_q^n$ and
$n-2$ zero columns.

For $d\in D$, let $b_1(d),\dots,b_n(d)$ be the columns of $X_0+d$, and
let
\[
 \varphi_d(s,t)=\det[\,s\ \ t\ \ b_3(d)\ \cdots\ b_n(d)\,],
\]
an alternating bilinear form on $\F_q^n$. Then
\[
\det(X_0+Y+d)=\varphi_d(y_1,y_2)+\varphi_d(y_1,b_2(d))
 +\varphi_d(b_1(d),y_2)+\varphi_d(b_1(d),b_2(d)).
\]
Expanding the product over $d\in D$, we obtain a sum of terms
\[
 g=c\prod_{d\in A}\varphi_d(y_1,y_2)\prod_{d\in B_1}\varphi_d(y_1,b_2(d))
 \prod_{d\in B_2}\varphi_d(b_1(d),y_2),\qquad c\in\F_q,
\]
one for each ordered partition $D=A\sqcup B_1\sqcup B_2\sqcup C$. Each
$g$ is a polynomial in the $2n$ coordinates of $y_1,y_2$, homogeneous of
degree $|A|+|B_1|$ in $y_1$ and of degree $|A|+|B_2|$ in $y_2$. We show
that the sum of every such $g$ over $y_1,y_2\in\F_q^n$ is zero.

Recall that $\sum_{x\in\F_q}x^e=0$ unless $e$ is a positive multiple of
$q-1$ (with $0^0=1$). Hence the sum of a monomial over $\F_q^{2n}$
vanishes unless each of its $2n$ exponents is a positive multiple of
$q-1$. Suppose that the sum of $g$ is nonzero. Then $|A|+|B_1|$ is a
sum of $n$ positive multiples of $q-1$, and
$|A|+|B_1|\le|D|<(n+1)(q-1)$. Thus $|A|+|B_1|=n(q-1)$, and likewise
$|A|+|B_2|=n(q-1)$. Put $b=|B_1|=|B_2|$. Then
$n(q-1)+b=|A|+2b\le|D|<(n+1)(q-1)$, so $b\le q-2$ and
$|A|=n(q-1)-b\ge1$.

Fix $y_1$. If $y_1=0$, then $g=0$ because $A\ne\varnothing$. If
$y_1\ne0$, split the sum over $y_2$ into the cosets $z+\F_qy_1$. Since
$\varphi_d$ is alternating, $\varphi_d(y_1,z+\lambda y_1)=\varphi_d(y_1,z)$,
so the factors indexed by $A$ and $B_1$ do not depend on $\lambda$,
while $\prod_{d\in B_2}\varphi_d(b_1(d),z+\lambda y_1)$ is a polynomial
in $\lambda$ of degree at most $b\le q-2$. The sum of such a polynomial
over $\lambda\in\F_q$ is zero. Hence the sum of $g$ over $y_2$ vanishes
for every $y_1$, a contradiction.

Finally, a summand $\prod_{d}\det(X_0+Y+d)$ is nonzero exactly when
$X_0+Y\in W(D)$. If there were exactly one such $Y$, the sum would be
nonzero.
\end{proof}

\begin{theorem}\label{Tm_n3}
Let $n\ge3$, let $q$ be a prime power, and let $D\subseteq M_n(\F_q)$ be a dominating set of
$\tau(M_n(\F_q))$ with $|D|\le n(q-1)$. Then $D$ is a singular-sum
covering set.
\end{theorem}

\begin{proof}
Put $W=W(D)$ and note that $W\subseteq D$. Suppose that $w_0\in W$. Since also
$w_0\in D$, the matrix $2w_0=w_0+w_0$ is invertible. If $q$ is even, then $2w_0=0$, a contradiction, so $W=\varnothing$. We may therefore from now on
assume that $q$ is odd.

\emph{Step 1.} By Lemma~\ref{L:sum}, applied with $X_0=w_0$ and any
$2$-dimensional subspace $R$, there is $w_1\in W$ with $w_1\ne w_0$ and
all rows of $w_1-w_0$ in $R$. Thus $\rank(w_1-w_0)\le2<n$, and
$w_1-w_0$ is singular.

\emph{Step 2.} Since $w_0,w_1\in W\subseteq D$, Lemma~\ref{L:coincide}(a)
with $d=w_0$, $d'=w_1$ and a nonzero vector $v$ in the kernel of
$w_1-w_0$ gives $r_1\ne0$ such that $x_1=w_0+u_1r_1^T\in W$, where
$u_1=2w_0v\ne0$.

\emph{Step 3.} Now $x_1\in D$, and the kernel of $x_1-w_0=u_1r_1^T$ is
the hyperplane $r_1^{\perp}=\{v:r_1^Tv=0\}$. Lemma~\ref{L:coincide}(a)
with $d=w_0$ and $d'=x_1$ shows that for every nonzero $u$ in the
$(n-1)$-dimensional space $V_1=2w_0\,r_1^{\perp}$ there is $r_u\ne0$ with
$x_u=w_0+ur_u^T\in W$.

\emph{Step 4.} Let $y\ne0$. Since $\dim V_1=n-1$ and $n\ge3$, the space
$V_1\cap y^{\perp}$ has dimension at least $n-2\ge1$. Choose $u\ne0$ in
it. Then $x_u\in D$ and $y^T(x_u-w_0)=(y^Tu)\,r_u^T=0$.
Lemma~\ref{L:coincide}(b) with $d=w_0$ and $d'=x_u$ gives $u'_y\ne0$
such that
\[
 z_y=w_0+u'_y\,(y^T2w_0)\in W .
\]
The matrix $z_y-w_0$ has rank one and row space spanned by $y^T2w_0$.
Since $2w_0$ is invertible, non-proportional vectors $y$ give
different row spaces, hence different matrices $z_y$, all distinct
from $w_0$.

\emph{Step 5.} Consequently,
\[
 |D|\ge|W|\ge1+\frac{q^n-1}{q-1}=1+\sum_{i=0}^{n-1}q^i
 \ge 1+n+\frac{n(n-1)}{2}(q-1)>n(q-1),
\]
where we used $q^i\ge1+i(q-1)$ and $n\ge3$. This contradicts
$|D|\le n(q-1)$. Hence $W=\varnothing$.
\end{proof}

The hypothesis $n\ge3$ is essential: for $n=2$ and odd $q$,
Proposition~\ref{Pr_upper} gives a dominating set of size
$2(q-1)$ that is not a covering set.

\begin{proof}[Proof of Theorem~\ref{Tm_n3_intro}]
The first assertion is Theorem~\ref{Tm_n3}.
Proposition~\ref{Prop_three}(b) gives
\[
  \gamma(\tau(M_n(\F_q)))\le\kappa(M_n(\F_q))\le n(q-1)+1.
\]
If the domination number is at most $n(q-1)$, then a minimum
dominating set is a covering set by Theorem~\ref{Tm_n3}, giving
the reverse inequality. Otherwise the domination number is
$n(q-1)+1$ and equality holds throughout. The final assertion
follows from~\eqref{For_conjecture}.
\end{proof}

For odd $q$ and $n\ge3$, Theorems~\ref{Tm_n3_intro}
and~\ref{Tm_shift} together give $\gamma=\kappa=\gamma_t$
for $\tau(M_n(\F_q))$.

The following corollary shows that, among the finite rings satisfying
the hypotheses of~\eqref{For_conjecture}, the domination and covering
numbers can differ only when $R$ itself is isomorphic to
$M_2(\F_q)$ for some odd $q$.

\begin{corollary}\label{Cor_products}
Let $R$ be a finite ring with
\[
 R/J(R)\cong\prod_{i=1}^k M_{n_i}(\F_{q_i}),\qquad n_i\ge2,
\]
and suppose that $R$ is not isomorphic to $M_2(\F_q)$ for any $q$. Then
\[
 \gamma(\tau(R))=\kappa(R)=\min_{1\le i\le k}\kappa(M_{n_i}(\F_{q_i})).
\]
Consequently, the value predicted by~\eqref{For_conjecture} is correct
for $R$ if and only if
$\min_i\kappa(M_{n_i}(\F_{q_i}))=\min_i\{n_i(q_i-1)+1\}$. In
particular, it is correct whenever the equality
$\kappa(M_{n_i}(\F_{q_i}))=n_i(q_i-1)+1$ holds for every factor.
\end{corollary}

\begin{proof}
Write $R_i=M_{n_i}(\F_{q_i})$. By Proposition~\ref{Pr_structure} and
induction on $k$, $\kappa(R)=\kappa(R/J(R))=\min_i\kappa(R_i)$. If
$J(R)\ne0$, the claim follows from Theorem~\ref{Tm_radical_domination},
and if $k=1$ and $J(R)=0$, then $R\cong R_1$ with $n_1\ge3$,
so the claim follows from Theorem~\ref{Tm_n3_intro}.
Assume therefore that $J(R)=0$ and $k\ge2$.
We identify $R$ with $R_1\times\cdots\times R_k$ via the given
isomorphism. By Proposition~\ref{Prop_three},
$\gamma(\tau(R))\le\kappa(R)$. 

Let $D\subseteq R$ be a dominating set with $|D|<\kappa(R)$. We derive
a contradiction. Let $D_i$ be the projection of $D$ to $R_i$. Since
$|D_i|\le|D|<\kappa(R_i)$, the set
\[
 W_i=\{x\in R_i: x+d\text{ is a unit for every }d\in D_i\}
\]
is nonempty. A sum in $R$ is a unit if and only if all its coordinates
are units, so $W_1\times\cdots\times W_k$ consists of vertices without a
neighbour in $D$, and therefore
\begin{equation}\label{For_product_W}
 W_1\times\cdots\times W_k\subseteq D .
\end{equation}
In particular $W_i\subseteq D_i$, that is, $D_i$ dominates $\tau(R_i)$
and $W(D_i)=W_i\ne\varnothing$, for every $i$. If some $n_i\ge3$, this
contradicts Theorem~\ref{Tm_n3}, since
$|D_i|\le\kappa(R_i)-1\le n_i(q_i-1)$. Hence $n_i=2$ for all $i$.

Since $\gamma(\tau(R_i))\le |D_i|<\kappa(R_i)$,
Proposition~\ref{Prop_three}(a) forces every $q_i$ to be odd.
Theorem~\ref{Tm_main}(b) now gives
$\gamma(\tau(R_1))=2q_1-2$. Thus $|D_1|\ge2q_1-2$, while $|D|\le\kappa(R_1)-1\le2q_1-2$ by
Proposition~\ref{Prop_three}(b). Hence $|D_1|=|D|$, and the projection
$D\to D_1$ is injective. On the other hand,
$|D_2|\le2q_2-2<3(q_2-1)$, so Lemma~\ref{L:sum}, applied with $n=2$, $D=D_2$, $X_0=0$
and the full row space $\F_{q_2}^2$, gives $|W_2|\ne1$.
Since $W_2$ is nonempty, we have $|W_2|\ge2$. Choose
$w_1\in W_1$, distinct $w_2,w_2'\in W_2$, and $w_i\in W_i$ for
$i\ge3$. By~\eqref{For_product_W}, the two distinct elements
$(w_1,w_2,w_3,\dots)$ and $(w_1,w_2',w_3,\dots)$ of $D$ have the same
first coordinate, contrary to injectivity. Hence
$\gamma(\tau(R))\ge\kappa(R)$.
\end{proof}

The hypothesis $n_i\ge2$ cannot be omitted. For
$R=\F_3\times\F_3$, the set $\{(0,0),(1,1)\}$ dominates
$\tau(R)$, whereas no single vertex dominates, since every
closed neighbourhood has at most $6$ elements. Thus
$\gamma(\tau(R))=2$, while Proposition~\ref{Pr_structure}
gives $\kappa(R)=\kappa(\F_3)=3$.

It remains to consider $R=M_2(\F_q)$. For even $q$,
Proposition~\ref{Prop_three}(a) gives $\gamma(\tau(R))=\kappa(R)$.
For odd $q$, Theorem~\ref{Tm_main}(b) gives
$\gamma(\tau(R))=\min\{\kappa(R),2q-2\}$, and
Proposition~\ref{Prop_three}(b) gives $\kappa(R)\le 2q-1$.
Hence a strict gap occurs exactly when $\kappa(R)=2q-1$,
in which case $\gamma(\tau(R))=\kappa(R)-1$.
By Theorem~\ref{Tm_main}(c), this occurs for $q=3,5,7$.

\bigskip

\noindent\textbf{Disclosure statement.}
The authors report there are no competing interests to declare.

\bigskip

\noindent\textbf{Declaration of AI use.}
The authors used Claude Claude Opus 5.5 (Anthropic) and ChatGPT 6.2 Astra (OpenAI) to help
develop and present mathematical arguments and to check proofs.
The authors verified all results and take full responsibility for
the content of the manuscript.

\end{document}